\documentclass[12pt]{amsart}
\usepackage[T1]{fontenc}
\usepackage[english]{babel}
\usepackage{color}
\usepackage{amssymb}
\usepackage{hyperref}
\usepackage{dsfont,soul}
\usepackage{enumitem}
\setlist[enumerate]{labelsep=*, leftmargin=1.5pc,
topsep=1ex plus0.5ex minus0.2ex,
itemsep=1ex plus0.5ex minus0.2ex,
font=\rmfamily,
font=\upshape}
\setlist[itemize]{labelsep=*, leftmargin=1.5pc,
topsep=1ex plus0.5ex minus0.2ex,
itemsep=1ex plus0.5ex minus0.2ex,
font=\rmfamily,
font=\upshape}
\newtheorem{thm}{Theorem}[section]
\newtheorem{cor}[thm]{Corollary}
\newtheorem{lem}[thm]{Lemma}
\newtheorem{pro}[thm]{Proposition}
\theoremstyle{definition}

\newtheorem{rem}[thm]{Remark}
\numberwithin{equation}{section}
\newcommand{\id}{\mathds{1}} 
\DeclareMathOperator{\ii}{i} 
\DeclareMathOperator{\tr}{tr}
\DeclareMathOperator{\diag}{diag}
\DeclareMathOperator*{\argmax}{argmax}

\newcommand{\m}{{\rm m}}
\newcommand{\bC}{\mathbb C}
\newcommand{\bE}{\mathbb E}
\newcommand{\bN}{\mathbb N}
\newcommand{\bP}{\mathbb P}
\newcommand{\bR}{\mathbb R}
\newcommand{\cE}{\mathcal E}
\newcommand{\cF}{\mathcal F}
\newcommand{\cM}{\mathcal M}
\newcommand{\cO}{\mathcal O}
\newcommand{\cR}{\mathcal R}
\begin{document}
\selectlanguage{english}
\author{Ilya M.\ Spitkovsky, Stephan Weis}
\title{Pre-images of extreme points of the numerical range,
and applications}
\begin{abstract}
We extend the pre-image representation of exposed points of the
numerical range of a matrix to all extreme points. With that we
characterize extreme points which are multiply generated, having
at least two linearly independent pre-images, as the extreme
points which are Hausdorff limits of flat boundary portions on
numerical ranges of a sequence converging to the given matrix.
These studies address the inverse numerical range map and the
maximum-entropy inference map which are continuous functions on
the numerical range except possibly at certain multiply generated
extreme points. This work also allows us to describe closures of
subsets of 3-by-3 matrices having the same shape of the numerical
range.
\end{abstract}
\keywords{Numerical range}
\subjclass[2010]{47A12, 54C10, 62F30, 94A17}
%
%
\maketitle
\thispagestyle{empty}
\pagestyle{myheadings}
\markleft{\hfill Pre-images of extreme points\hfill}
\markright{\hfill I.\,M.~Spitkovsky, S.~Weis \hfill}
%
%
\section{Introduction}
\par
We denote the set of complex $d\times d$ matrices by $M_d$,
$d\in\bN$, with identity matrix $\id$. The
{\em numerical range} of $A\in M_d$ is the subset
\[
W(A):=\{f_A(x)\mid x\in S\bC^d\}
\]
of the complex plane $\bC$ where $f_A:S\bC^d\to\bC,x\mapsto x^*Ax$
is the {\em numerical range map} defined on the unit sphere
$S\bC^d=\{x\in\bC^d\mid x^*x=1\}$ of $\bC^d$. On $\bC\cong\bR^2$
we use the standard Euclidean scalar product
$\langle\alpha,\beta\rangle=\Re(\alpha\overline{\beta})$ for
$\alpha,\beta\in\bC$. The numerical range is a compact and convex
subset, the convexity statement being known as the Toeplitz-Hausdorff
theorem \cite{Toeplitz1918,Hausdorff1919}, more recent work
includes \cite{Au-YeungPoon1979,HornJohnson1991,LiPoon2000}.
\par
We are interested in points $\alpha\in W(A)$ which are
\begin{itemize}
\item
{\em extreme points}, that is $\alpha$ cannot be written as a proper
convex combination of points in $W(A)$,
\item
{\em multiply generated}
\cite{Leake-etal-OAM2014,Leake-etal-LMA2014,Rodman-etal2016}, that is
$f_A^{-1}(\alpha)$ contains at least two linearly independent vectors.
\end{itemize}
\par
Part of our interest in these points comes from a continuity problem
in operator theory. The (multi-valued) inverse $f_A^{-1}: W(A)\to S\bC^d$
of $f_A$ is {\em strongly continuous} at $\alpha\in W(A)$ if for all
$x\in f_A^{-1}(\alpha)$ the function $f_A$ is open at $x$. Strong
continuity holds on $W(A)$ except at certain multiply generated extreme
points. A {\em round boundary point}
(see \cite{Corey-etal2013,Leake-etal-OAM2014,Leake-etal-LMA2014} and
Lemma~6.1 of \cite{Rodman-etal2016})
is an extreme point of $W(A)$ which lies on a unique {\em supporting line}
$\{\alpha\in\bC\mid\langle\alpha,e^{\ii\theta}\rangle=h(e^{\ii\theta})\}$
where $e^{\ii\theta}$, $\theta\in\bR$, is an outward pointing normal vector
and
$h(e^{\ii\theta})=\max_{\alpha\in W(A)}\langle\alpha,e^{\ii\theta}\rangle$.
The map $f_A^{-1}$ is strongly continuous on $W(A)$ except possibly at
multiply generated round boundary points \cite{Corey-etal2013} but may
be strongly continuous also there. A characterization of strong continuity
of $f_A^{-1}$, in terms of analytic eigenvalue curves
\cite{Leake-etal-LMA2014} shows that $f_A^{-1}$ has at most finitely many
discontinuity points. The corresponding eigenvector curves will be discussed
in Coro.~\ref{cor:ev-curves-rep} where we consider their intersection with
pre-images of extreme points under $f_A$.
\par
Further interest in multiply generated extreme points comes from an optimization problem
in quantum mechanics. The
{\em maximum-entropy inference}, going back to ideas by Boltzmann, is a
method to select a quantum state from the expected values of a collection
of physical quantities represented by hermitian matrices when all other
information about the state is ignored \cite{Jaynes1957}. The numerical
range $W(A)$ is the set of expected values of two hermitian
matrices given implicitly by the {\em real part}
$\Re(A)=\tfrac{1}{2}(A+A^*)$ and {\em imaginary part}
$\Im(A)=\tfrac{1}{2\ii}(A-A^*)$ of $A$
\cite{BerberianOrland1967,Rodman-etal2016}. The inference map is
\[
\rho^*:W(A)\to\cM_d,\quad
\alpha\mapsto\argmax\{S(\rho)\mid\rho\in\cM_d,\tr(A\rho)=\alpha\}
\]
with {\em state space} $\cM_d$ consisting of all positive
semi-definite matrices of trace one and {\em von Neumann entropy}
$S(\rho)=-\tr\rho\log(\rho)$. A discontinuity of $\rho^*$ may occur
\cite{WeisKnauf2012}. Note, however, that for normal $A$ its
numerical range $W(A)$ is a polytope, $\Re(A)$ and $\Im(A)$ commute,
and $\rho^*$ is continuous \cite{Weis-cont}. So, a discontinuity
belongs to the proper quantum domain where it is discussed in the
context of quantum phases \cite{Chen-etal2015,Kato-etal2016}.
\par
As it happens, the
described continuity problems are equivalent. The map $\rho^*$ is
continuous at $\alpha\in W(A)$ if and only if $f_A^{-1}$ is strongly
continuous at $\alpha$ \cite{Weis-strong}, and $\rho^*$ is indeed
discontinuous at all {\em isolated} multiply generated round boundary
points of $W(A)$, see Sec.~6 of \cite{Rodman-etal2016}. Further, it is
known that multiply generated round boundary points are isolated for
$d=3$. Calling $A\in M_d$ {\em unitarily reducible}
if $A$ is unitarily similar to a block diagonal matrix with two proper
blocks and {\em unitarily irreducible} otherwise,
multiply generated round boundary points are isolated for all
irreducible matrices of size $d\leq 5$ \cite{Leake-etal-OAM2014}. But,
the whole boundary of $W(A)$ may consist of multiply generated round
boundary points for reducible 4-by-4 matrices and irreducible 6-by-6
matrices \cite{Leake-etal-OAM2014}.
\par
Here we study these continuity problems from a topological perspective
of a variable matrix $A\in M_d$. To this end we define a
{\em flat boundary portion} of $W(A)$ as
a maximal proper segment in the boundary of $W(A)$
\cite{Keeler-etal1997,BrownSpitkovsky2004,RodmanSpitkovsky2005,Eldred-etal2012}.
If $\alpha\in W(A)$, then we say that
{\em a flat boundary portion is born at $\alpha$} if there exists a
sequence $(A_i)_{i\in\bN}\subset M_d$ converging to $A$ such that a
sequence $(s_i)_{i\in\bN}$ of flat boundary portions converges to
$\{\alpha\}$ in the Hausdorff distance, the segment $s_i$ being a flat
boundary portion of $W(A_i)$, $i\in\bN$.
\par
The birth of a flat boundary portion was conjectured in Sec.~I.B of
\cite{WeisKnauf2012} to be a condition for a discontinuity of
$\rho^*$. The above discussion and the following theorem prove that
the birth of a flat boundary portion at a round boundary point is a
necessary condition. It is not a sufficient condition for $d\geq 4$ 
($d\geq 6$ if $A$ is unitarily irreducible).
\begin{thm}\label{thm:birth}
Let $\alpha$ be an extreme point of $W(A)$.
For a flat boundary portion to be born at $\alpha$ it is necessary
and sufficient that $\alpha$ is multiply generated.
\end{thm}
The necessity in Thm.~\ref{thm:birth} follows from properties of the
Hausdorff distance (Sec.~\ref{sec:Hausdorff}). We prove the sufficiency
in Sec.~\ref{sec:eigen-rep} using a newly developed representation of
extreme points in terms of pre-images. Pre-images were well-understood
\cite{Toeplitz1918} for {\em exposed points} which can be represented
as the intersection of $W(A)$ with a supporting line. The number
$h(e^{\ii\theta})$ defined above is the maximal eigenvalue of the
hermitian matrix $\cos(\theta)\Re(A)+\sin(\theta)\Im(A)$, the
corresponding eigenspace is the demanded pre-image.
\par
Some extreme points may fail to be exposed. Consider, e.g., the convex
hull of a circle and a point outside the circle,
which is realized as the numerical range of
$A=\left[\begin{smallmatrix}
0 & 2 & 0\\
0 & 0 & 0\\
0 & 0 & 2
\end{smallmatrix}\right]$.
Then the intersections of the two tangents from the point to the circle
with the circle are
extreme but non-exposed points. We will obtain the pre-images of
non-exposed points by viewing them as exposed points of some flat
boundary portion whose supporting lines are given by directional
derivatives of $h$ \cite{BonnesenFenchel1987}. Notice that viewing
non-exposed points as exposed points of a convex subset is a familiar
idea in convex geometry
\cite{Gruenbaum2003,Weis-cones} (cf.~{\em poonem}), geometry of
quantum states \cite{Weis-support}, and in the theory of exponential
families \cite{CsiszarMatus2005,Weis-topo} (cf.~{\em access sequence}).
\par
Finally, in Sec.~\ref{sec:3x3matrices} we combine observations from
\cite{Halmos1970,Kippenhahn1951,Keeler-etal1997,RodmanSpitkovsky2005,%
Leake-etal-OAM2014,Rault-etal2013} with Thm.~\ref{thm:birth} to compute
closures of subsets of 3-by-3 matrices having the same shape of the
numerical range. In particular, we prove in Sec.~\ref{sec:reducible}
that a 3-by-3 matrix $A$ with $W(A)$ having a non-empty interior lies in
the closure of irreducible matrices with elliptical numerical range and
in the closure of irreducible matrices with flat portion on the boundary
of their numerical range if and only if $W(A)$ is an ellipse with an
eigenvalue of $A$ on the
boundary. Remarkably, these are precisely the 3-by-3 matrices $A$
where $f_A^{-1}$ and $\rho^*$ have a discontinuity
\cite{Corey-etal2013,Leake-etal-OAM2014,Rodman-etal2016}.\\
\par\normalsize
{\par\noindent\footnotesize
{\em Acknowledgements.}
SW cherishes the memories of discussions with Leiba Rodman during
the workshop LAW'14 in Ljub\-lja\-na, Slovenia, June 4--12, 2014.
He also appreciates financial support by a Brazilian Capes scholarship.
IS acknowledges the support by the Plumeri Award for
Faculty Excellence from the College of William and Mary and by Faculty
Research funding from the Division of Science and Mathematics, New York
University Abu Dhabi (NYUAD). Both authors are especially thankful to
the latter for hosting SW's research visit to NYUAD in December'14.}
%
%
%
\section{Representation of extreme points}
\label{sec:eigen-rep}
\par
We provide a pre-image representation of extreme points of the
numerical range. We rewrite it in terms of eigenvalue curves and
we use it to prove one part of Thm.~\ref{thm:birth}.
\par
Recall that the {\em relative interior} of a subset $C$ of
a Euclidean space $(\bE,\langle\cdot,\cdot\rangle)$ is the interior
of $C$ in the affine hull of $C$. A {\em face} of a convex set
$C$ is a convex subset $F\subset C$ which contains every closed
segment in $C$ the relative interior of which intersects $F$.
The element of a singleton face is called {\em extreme point}.
The {\em support function} of a {\em convex body} $C$, that is a
compact convex subset of $\bE$, is
\[
h_C(u):=\max_{x\in C}\langle x,u\rangle,
\qquad u\in\bE.
\]
The {\em supporting hyperplane} of $C$ with outward pointing normal
vector $u\in\bE\setminus\{0\}$ is the set of $x\in\bE$ on the hyperplane
$\langle x,u\rangle=h_C(u)$. The intersection of this hyperplane with
$C$ is the {\em exposed face}
\[
\argmax_{x\in C}\langle x,u\rangle
\]
with outward pointing normal vector $u$. If $\{\alpha\}$ is an exposed
face, then $\alpha$ is called an {\em exposed point}. All exposed points
are extreme points. The remaining extreme points are called
{\em non-exposed points}.
\par
Notice that a supporting hyperplane of $W(A)$ is a supporting line. Since
$W(A)$ has (real) dimension at most two, a subset is a one-dimensional
face if and only if it is a flat boundary portion, and every flat 
boundary portion is an exposed face. Moreover, the boundary of $W(A)$ is 
the disjoint union of extreme points and relative interiors of flat 
boundary portions. So, every extreme point $\alpha$ is an exposed point, 
an endpoint of a flat boundary portion, or both. It is crucial in the 
sequel that every non-exposed point of $W(A)$ is an exposed point of some 
flat boundary portion of $W(A)$. 
\par
Let
\[
A(\theta):=\Re(e^{-\ii\theta}A)=\cos(\theta)\Re(A)+\sin(\theta)\Im(A),
\qquad
\theta\in\bR.
\]
We denote by $X_\m(\theta)$ the eigenspace of $A(\theta)$
corresponding to the maximal eigenvalue $\lambda_\m(\theta)$
of $A(\theta)$. An easy computation gives
\begin{equation}\label{eq:rad-coord}
\langle f_A(x),e^{\ii\theta}\rangle
=\Re(x^*Axe^{-\ii\theta})
=f_{A(\theta)}(x), \qquad x\in S\bC^d.
\end{equation}
The maximal eigenvalue $\lambda_\m(\theta)$ of $A(\theta)$ has the
geometric meaning of support function of the numerical range, see
Sec.~4 of \cite{Toeplitz1918},
\begin{equation}\label{eq:support-W}
h_{W(A)}(e^{\ii\theta})=\lambda_\m(\theta).
\end{equation}
Let $F_A(\theta)$ denote the exposed face of $W(A)$
with outward pointing normal vector $e^{\ii\theta}$, $\theta\in\bR$.
\begin{lem}\label{lem:exp-face-rep}
The point $f_A(x)$, $x\in S\bC^d$, lies in the
exposed face $F_A(\theta)$, $\theta\in\bR$, if and only if
$x\in X_\m(\theta)$.
\end{lem}
{\em Proof:}
We consider the orthogonal direct sum
$\bC^d=X_\m(\theta)\oplus X_\m(\theta)^\perp$. A unit vector
$x\in S\bC^d$ has a unique decomposition $x=y+z$ for
$y\in X_\m(\theta)$ and $z\in X_\m(\theta)^\perp$. So
\[
\langle f_A(x),e^{\ii\theta}\rangle
=\lambda_\m(\theta)+z^*(A(\theta)-\lambda_\m(\theta)\id)z
\]
has the maximal value $\lambda_\m(\theta)$ only for $z=0$. This
proves the claim.
\hspace*{\fill}$\square$\\
\par
Non-exposed points of the numerical range are not addressed in
Lemma~\ref{lem:exp-face-rep}. To describe them in terms of
pre-images, we view them as exposed points of a flat boundary
portion. The two endpoints $p_\pm(\theta)$ (not necessarily distinct)
of the exposed face $F_A(\theta)$ are characterized by their
membership in $F_A(\theta)$ and by the equality
\[
\langle p_\pm(\theta),\pm\ii e^{\ii\theta}\rangle
=h_{F_A(\theta)}(\pm\ii e^{\ii\theta}).
\]
To evaluate the support function of $F_A(\theta)$ we consider
the directional derivative of a function $f:\bR^m\to\bR^n$ at a
point $u\in\bR^m$ in direction $v\in\bR^m$, which is defined by
\[
f'(u;v):=\lim_{t\to0^+}\tfrac{1}{t}(f(u+tv)-f(u)),
\]
if the limit exists. Given a convex body $C$ in a Euclidean
space $\bE$ and $u\in\bE\setminus\{0\}$, the support function
$h_G$ of the exposed face $G$ of $C$ with outward pointing normal
vector $u$ is $h_G(v)=h_C'(u;v)$, $v\in\bE$, see Section~16
of \cite{BonnesenFenchel1987}.
\par
Since the chain rule does not hold for directional derivatives,
we need to go into some detail in the following proof.
\begin{lem}\label{lem:face-rep-directional}
The support function of the exposed face $F_A(\theta)$, $\theta\in\bR$,
has the values
$h_{F_A(\theta)}(\pm\ii e^{\ii\theta})=\lambda_\m'(\theta;\pm 1)$.
The endpoints of $F_A(\theta)$ are
$p_\pm(\theta)=e^{\ii\theta}(\lambda_\m(\theta)\pm\ii\lambda_\m'(\theta;\pm 1))$.
\end{lem}
{\em Proof:}
The equation for support functions of exposed faces, recalled in the previous
paragraph, proves
\[
h_{F_A(\theta)}(\pm\ii e^{\ii\theta})
=h_{W(A)}'(e^{i\theta};\pm\ii e^{\ii\theta}).
\]
For all $s,t>0$ with $s=\arctan(t)$, and writing $h=h_{W(A)}$, we have
\[
\tfrac{1}{t}(h(e^{\ii\theta}\pm t\ii e^{\ii\theta})-h(e^{\ii\theta}))
=\tfrac{1}{s}(h(e^{\ii(\theta\pm s)})-h(e^{\ii\theta}))+g(s)
\]
where $g(s)\to 0$ for $s\to 0$, because $h$ is positive homogeneous of
the first degree \cite{BonnesenFenchel1987}.
Taking the limit $t\to 0^+$ gives
\[
h_{W(A)}'(e^{\ii\theta};\pm\ii e^{\ii\theta})
=(h_{W(A)}\circ e^{\ii\theta})'(\theta;\pm 1).
\]
The values of $h_{F_A(\theta)}$ follow from
\eqref{eq:support-W} which provides
$h_{W(A)}\circ e^{\ii\theta}=\lambda_\m(\theta)$.
The formula for the endpoints is obvious.
\hspace*{\fill}$\square$\\
\par
We are ready to describe $p_\pm(\theta)$ in terms of its pre-image
under $f_A$. Noticing $A'(\theta)=\Im(e^{-\ii\theta}A)$, where
``\,$'$\,'' denotes derivative with respect to $\theta$, an easy
calculation shows that for all $y\in S\bC^d$
\begin{equation}\label{eq:coordinate-Aprime}
\langle f_A(y),\ii e^{\ii\theta}\rangle=f_{A'(\theta)}(y)
\end{equation}
holds. Equations
\eqref{eq:rad-coord} and \eqref{eq:coordinate-Aprime} show that for
all $y\in S\bC^d$ we have
\begin{equation}\label{eq:fAofy}
f_A(y)=e^{\ii\theta}(f_{A(\theta)}(y)+\ii f_{A'(\theta)}(y)).
\end{equation}
We denote by $B|_X$ the compression of $B\in M_d$ onto a subspace
$X\subset\bC^d$, that is $B|_X$ is the restriction of $PBP$ to $X$
where $P$ is the orthogonal projection onto $X$. For $\theta\in\bR$ define
$g_\theta:\bC\to\bR,\alpha\mapsto\langle \alpha,\ii e^{\ii\theta}\rangle$
and
$h_\theta:\bR\to\bC,\eta\mapsto e^{\ii\theta}(\lambda_\m(\theta)+\ii\eta)$.
\begin{thm}\label{eq:thm-face-rep}
For all $\theta\in\bR$ the maximal (respectively, minimal) eigenvalue
of $A'(\theta)|_{X_\m(\theta)}$ is $\lambda_\m'(\theta;1)$
(respectively, $-\lambda_\m'(\theta;-1)$). For all unit vectors
$x\in X_\m(\theta)$ we have
$f_{A'(\theta)|_{X_\m(\theta)}}(x)=g_\theta\circ f_A(x)$. The map
$g_\theta|_{F_A(\theta)}:F_A(\theta)\to s$ is a bijection to the
segment $s\subset\bR$ with endpoints $\pm\lambda_\m'(\theta;\pm1)$,
the inverse is $h_\theta|_s$.
\end{thm}
{\em Proof:}
The pre-image of the exposed face $F_A(\theta)$ is by
Lemma~\ref{lem:exp-face-rep} equal to $S\bC^d\cap X_\m(\theta)$.
By definition of the support function of $F_A(\theta)$ we get for
all $x\in S\bC^d\cap X_\m(\theta)$ the inequalities
\[
-h_{F_A(\theta)}(-\ii e^{\ii\theta})
\leq\langle f_A(x),\ii e^{\ii\theta}\rangle
\leq h_{F_A(\theta)}(\ii e^{\ii\theta}).
\]
Both equalities are attained because $F_A(\theta)$ is compact. Using
Lemma~\ref{lem:face-rep-directional} and
\eqref{eq:coordinate-Aprime}, the above inequality is equivalent to
\[
-\lambda_\m'(\theta;-1)
\leq x^*A'(\theta)x
\leq\lambda_\m'(\theta;+1)
\]
which shows that the hermitian operator $A'(\theta)|_{X_\m}$ has
minimal eigenvalue $-\lambda_\m'(\theta;-1)$ and maximal
eigenvalue $\lambda_\m'(\theta;+1)$.
\par
The numerical range of the hermitian operator
$A'(\theta)|_{X_\m(\theta)}$ is the segment $s$ between its
extreme eigenvalues $\pm\lambda_\m'(\theta;\pm1)$. By
\eqref{eq:support-W} and \eqref{eq:fAofy} we have for all
$x\in S\bC^d\cap X_\m(\theta)$
\[
f_A(x)=e^{\ii\theta}(\lambda_\m(\theta)+\ii f_{A'(\theta)}(x)).
\]
Since $f_{A'(\theta)}(x)=\langle f_A(x),\ii e^{\ii\theta}\rangle$
holds, again by \eqref{eq:coordinate-Aprime}, the functions
$g_\theta$ and $h_\theta$ have the claimed properties.
\hspace*{\fill}$\square$\\
\par
We can now compute the pre-images of all extreme points.
Thm.~\ref{eq:thm-face-rep} proves that
$\pm\lambda_\m'(\theta;\pm1)$ is an eigenvalue of $A'(\theta)|_{X_\m}$.
We denote the corresponding eigenspace by $X_\pm(\theta)$.
\begin{cor}\label{cor:pre-image-p}
The point $f_A(x)$, $x\in S\bC^d$, is the endpoint $p_\pm(\theta)$,
$\theta\in\bR$, of the exposed face $F_A(\theta)$ if and only if
$x\in X_\pm(\theta)$.
\end{cor}
{\em Proof:}
This follows from Thm.~\ref{eq:thm-face-rep}, and from
Lemma~\ref{lem:exp-face-rep} applied to the exposed points
$\pm\lambda'_\m(\theta;\pm1)$ of the numerical range of the
compression $A'(\theta)|_{X_\m(\theta)}$.
\hspace*{\fill}$\square$\\
\par
We describe the eigenspace $X_\pm(\theta)$ of $A'(\theta)|_{X_\m}$
in terms of eigenvalue curves $\{\lambda_k(\theta)\}_{k=1}^d$ of
$A(\theta)$ and mutually orthogonal eigenvectors $\{x_k(\theta)\}_{k=1}^d$ which
depend real analytically on the parameter $\theta$
\cite{Rellich1954}. Thus
$A(\theta)=\sum_{k=1}^d\lambda_k(\theta)x_k(\theta)x_k(\theta)^*$
and we have for $k=1,\ldots,d$
\[
\lambda_k(\theta)=f_{A(\theta)}(x_k(\theta)),
\qquad\theta\in\bR.
\]
We recall that the derivative of $\lambda_k$ with respect to $\theta$
(see Lemma~3.2 of \cite{JoswigStraub1998} and Sec.~5 of \cite{GallaySerre2012})
is
\begin{equation}\label{eq:JoswigStraub}
\lambda_k'(\theta)=f_{A'(\theta)}(x_k(\theta)),
\qquad\theta\in\bR.
\end{equation}
We define
\[
I_\m(\theta)
:=\{i\in\{1,\ldots,d\}\mid\lambda_i(\theta)=\lambda_\m(\theta)\}
\]
and
\[
I_\pm(\theta)
:=\{i\in I_\m(\theta)\mid
\lambda'_i(\theta)=\pm\lambda_\m'(\theta;\pm1)\}.
\]
\begin{cor}\label{cor:ev-curves-rep}
The span of $\{x_k(\theta)\mid k\in I_\pm(\theta)\}$
is $X_\pm(\theta)$, $\theta\in\bR$.
\end{cor}
{\em Proof:}
By Thm.~\ref{eq:thm-face-rep} the eigenvalue $\pm\lambda_\m'(\theta;\pm1)$
of the hermitian operator $A'(\theta)|_{X_\m(\theta)}$ is extreme.
Hence, \eqref{eq:JoswigStraub} shows that $X_\pm(\theta)$ contains
all eigenvectors $x_k(\theta)$ of $A(\theta)$ with
$k\in I_\pm(\theta)$.
For all $k\in I_\m(\theta)\setminus I_\pm(\theta)$ we have by
\eqref{eq:JoswigStraub}
\[
x_k(\theta)^*A'(\theta)x_k(\theta)
=\lambda_k'(\theta)
\lessgtr\pm\lambda_\m'(\theta;\pm1).
\]
Therefore $X_\pm(\theta)$ is the span of
$\{x_k(\theta)\mid k\in I_\pm(\theta)\}$.
\hspace*{\fill}$\square$\\
\par
We finish the section with an application. See the next
section for the definition of the Hausdorff distance.
\begin{pro}\label{pro:birth-of-fbp}
Let $\alpha$ be a multiply generated extreme point of the numerical
range $W(A)$. Then a flat boundary portion is born at $\alpha$.
\end{pro}
{\em Proof:}
Without loss of generality let $\alpha=p_\sigma(\theta)$ for some
$\theta\in\bR$ and some $\sigma\in\{+,-\}$. Since $X_\sigma(\theta)$
belongs to the eigenspace $X_\m(\theta)$ of $A(\theta)$ corresponding
to the maximal eigenvalue $\lambda_\m(\theta)$ of $A(\theta)$, the
equation \eqref{eq:fAofy} implies that for all
$x\in S\bC^d\cap X_\sigma(\theta)$
\[
f_A(x)=e^{\ii\theta}(\lambda_\m(\theta)+\ii f_{A'(\theta)}(x)).
\]
Let $P$ denote the orthogonal projection onto $X_\sigma(\theta)$.
Coro.~\ref{cor:pre-image-p} shows that $X_\sigma(\theta)$ is the
pre-image of $\alpha$ under $f_A$. Hence, there exists $\lambda\in\bR$
such that $A'(\theta)|_{X_\sigma(\theta)}=\lambda P|_{X_\sigma(\theta)}$,
for otherwise $f_A(X_\sigma(\theta))$ could not be a singleton. As
the numerical range of $\lambda P|_{X_\sigma(\theta)}$ is $\{\lambda\}$,
Thm.~\ref{eq:thm-face-rep} proves
$\alpha=e^{\ii\theta}(\lambda_\m(\theta)+\ii \lambda)$.
\par
By assumption, $\alpha$ is multiply generated, so
$\dim_\bC(X_\sigma(\theta))\geq 2$ holds. Choose any hermitian matrix $H$
with $H|_{X_\sigma(\theta)}$ not being a scalar multiple of the identity
and denote its maximal eigenvalue by $\mu_+$ and minimal eigenvalue by
$\mu_-$. Let $\epsilon>0$, and define
\[
A_\epsilon:=e^{\ii\theta}
(A(\theta)+\epsilon P+\ii(A'(\theta)+\epsilon H)).
\]
The numerical range of $(A'(\theta)+\epsilon H)|_{X_\sigma(\theta)}$ is
the segment between the two distinct reals $\lambda+\epsilon\mu_\pm$. Since
$X_\sigma(\theta)$ is the eigenspace of $\Re(e^{-\ii\theta}A_\epsilon)$
corresponding to the maximal eigenvalue $\lambda_\m(\theta)+\epsilon$ of
$\Re(e^{-\ii\theta}A_\epsilon)$, the exposed face
$F_{W(A_\epsilon)}(\theta)$ of $W(A_\epsilon)$ has by
Thm.~\ref{eq:thm-face-rep} the endpoints
\[
e^{\ii\theta}(\lambda_\m(\theta)+\epsilon+\ii(\lambda+\epsilon\mu_\pm)).
\]
Clearly, $F_{W(A_\epsilon)}(\theta)$ is a flat boundary portion of
$W(A_\epsilon)$ which converges for $\epsilon\to 0$ in the Hausdorff
distance to $\{\alpha\}$ while $A_\epsilon$ converges to $A$. This
completes the proof.
\hspace*{\fill}$\square$\\
%
%
\section{Hausdorff distance}
\label{sec:Hausdorff}
\par
We address the Hausdorff convergence of numerical ranges and
show that a flat boundary portion can only be born at a point
of the numerical range which is multiply generated.
\par
We denote by $|\cdot|$ the Euclidean norm in a Euclidean space
$(\bE,\langle\cdot,\cdot\rangle)$. The set of non-empty compact
subsets of $\bE$ is a complete metric space with respect to the
{\em Hausdorff distance}
\[
d_H(K,L):=
\max(\max_{x\in K}\min_{y\in L}|y-x|,
\max_{y\in L}\min_{x\in K}|y-x|),
\]
where $\emptyset\neq K,L\subset\bE$ are compact. The set of all
convex bodies in $\bE$ is closed in the metric space of compacts
by Thm.~1.8.5 in \cite{Schneider2014}.
\par
We first recall a condition for the Hausdorff convergence.
\begin{rem}[Thm.~1.8.7 in \cite{Schneider2014}]\label{rem:Schneider}
A sequence $(K_i)_{i\in\bN}$ of convex bodies in $\bE$ converges
to a convex body $K$ in $\bE$ if and only if the following two
conditions hold.
\begin{enumerate}
\item Each point in $K$ is the limit of a sequence $(x_i)_{i\in\bN}$
with $x_i\in K_i$ for $i\in\bN$.
\item
The limit of any convergent subsequence $(x_{i_j})_{j\in\bN}$ with
$x_{i_j}\in K_{i_j}$ for $j\in\bN$ belongs to $K$.
\end{enumerate}
\end{rem}
The simplest example of a non-converging sequence illustrates that
subsequences are needed in (2). Consider for example $K=\{0\}$ and
the sequence defined by $K_{2i}=K$ and $K_{2i+1}=[0,1]$, $i\in\bN$.
\begin{rem}[Hausdorff convergence of the numerical range]%
\label{rem:Hcnr}~
\begin{enumerate}
\item
Remark 3.1 and the continuity for every $x\in\bC^d$ of the linear
function $M_d\to\bC$, $A\mapsto x^*Ax$ in $A$ prove that
\[
A_i\stackrel{i\to\infty}{\longrightarrow}A
\quad\implies\quad
W(A_i)\stackrel{i\to\infty}{\longrightarrow}W(A),
\qquad (A_i)_{i\in\bN}\subset M_d.
\]
Here, the first limit is in any norm, the second limit is in the
Hausdorff distance.
\item
The support function \eqref{eq:support-W} of $W(A)$ in the direction
$e^{\ii\theta}$ is jointly continuous in $A\in M_d$ and $\theta\in\bR$
because the maximal eigenvalue of a hermitian matrix is a continuous
function of the matrix. Therefore, if $\theta=\lim_{i\to\infty}\theta_i$,
$(\theta_i)_{i\in\bN}\subset\bR$, $A=\lim_{i\to\infty}A_i$,
$(A_i)_{i\in\bN}\subset M_d$, and if the exposed faces $F_{A_i}(\theta_i)$
converge in the Hausdorff distance, then their limit is a subset of
$F_A(\theta)$.
\end{enumerate}
\end{rem}
\par
We now complete the proof of Thm.~\ref{thm:birth}.
Recall from Section~\ref{sec:eigen-rep}, third paragraph,
that a flat boundary portion of $W(A)$ is an exposed face 
of dimension one.
\begin{lem}\label{lem:birth-implies-multiple}
If a flat boundary portion is born at $\alpha\in W(A)$, then
$\alpha$ is multiply generated.
\end{lem}
{\em Proof:} Let $(s_i)_{i\in\bN}$ be a sequence of segments,
each $s_i$ being a flat boundary portion of $W(A_i)$, $i\in\bN$,
and assume that $A=\lim_{i\to\infty}A_i$ and that
$s_i\stackrel{i\to\infty}{\rightarrow}\{\alpha\}$
converges in the Hausdorff distance. The segments $s_i$ are
exposed faces of $W(A_i)$. So, for every $i\in\bN$ there exists
$\theta_i\in\bR$ such that $e^{\ii\theta_i}$ is an outward pointing
normal vector of $s_i$. Lemma~\ref{lem:exp-face-rep} implies
\[
f_{A_i}^{-1}(s_i)=S\bC^d\cap H_i
\]
for some subspace $H_i\subset\bC^d$. The segments $s_i$ are no
singletons so $\dim_\bC(H_i)\geq 2$. This proves that
$f_{A_i}^{-1}(s_i)$ contains two orthogonal vectors $p_i$,
$q_i$. The compactness of $S\bC^d$ proves
$p_{i_j}\stackrel{j\to\infty}{\rightarrow}p$ and
$q_{i_j}\stackrel{j\to\infty}{\rightarrow}q$ for suitable 
$p,q\in S\bC^d$ and a suitable subsequence. Note that $p$ and $q$ 
are orthogonal. Since $f_A(p)$ is jointly continuous in $A$ 
and $p$, and since $f_{A_i}(p_i)\in s_i$ for all $i\in\bN$, this shows
\begin{align*}
|f_A(p)-\alpha| 
 & =\lim_{j\to\infty}|f_{A_{i_j}}(p_{i_j})-\alpha|
 \leq \limsup_{j\to\infty}\max_{\beta\in s_{i_j}}|\beta-\alpha|\\
 & \leq \limsup_{j\to\infty}d_H(s_{i_j},\{\alpha\})=0.
\end{align*}
Similarly, $f_A(q)=\alpha$ holds and hence $\alpha$ is multiply
generated.
\hspace*{\fill}$\square$\\
\par
Before turning to other subjects we add a statement about the Hausdorff
convergence of ellipses in the plane. Thereby an {\em ellipse} is the
zero set of a real quadratic form $q:\bR^2\to\bR$,
$x\mapsto x^*Sx+b^*x+c$ where $S\in M_2$ is a real symmetric and
positive definite 2-by-2 matrix, $b\in\bR^2$, and $c\in\bR$. We
call the level set $\{x\in\bR^2\mid f(x)\leq 0\}$, which is the convex
hull of the ellipse, also an ellipse.
\begin{rem}\label{rem:limit-of-ellipses}
If a sequence of ellipses in $\bR^2$ converges in the Hausdorff distance
then the limit is an ellipse, a segment, or a point. This can be
proved by representing each ellipse as a linear image of the unit disk
and by using Rem.~\ref{rem:Schneider} and compactness arguments.
\end{rem}
%
%
\section{3-by-3 matrices}
\label{sec:3x3matrices}
\par
Kippenhahn's \cite{Kippenhahn1951} representation of the numerical range
of a $d\times d$ matrix $A$ in terms of the convex hull of a planar real
affine algebraic curve provides a classification of the possible shapes
of the numerical range of a 3-by-3 matrix whose equivalence classes
are well-understood \cite{Keeler-etal1997,RodmanSpitkovsky2005} in terms
of matrix entries and of spectral data of $A$. We compute the closures of
these equivalence classes.
\par
Using the hermitian real and imaginary parts of $A=\Re(A)+\ii\Im(A)$,
the homogeneous polynomial
\[
F(y_0,y_1,y_2):=\det(y_0\id+y_1\Re(A)+y_2\Im(A))
\]
defines a curve $S_A:=\{(y_0:y_1:y_2)\in\bP\bC^3\mid F(y_0,y_1,y_2)=0\}$
in the complex projective plane $\bP\bC^3$. The dual curve $S_A^\wedge$
is an algebraic curve in $\bP\bC^3$ which consists roughly speaking of
the tangent lines to $S_A$ \cite{Gelfand-etal1994,Fischer2001}.
Thereby, a line $\{(y_0:y_1:y_2)\in\bP\bC^3\mid x_0y_0+x_1y_1+x_2y_2=0\}$
is identified with the point $(x_0:x_1:x_2)\in\bP\bC^3$. The
{\em boundary generating curve} $S_A^\wedge(\bR)$ of $A$ is the real part
of the affine component $x_0=1$ of $S_A^\wedge$. One can show that the 
convex hull of $S_A^\wedge(\bR)$ is the numerical range
\cite{ChienNakazato2010,Kippenhahn1951}. See also
\cite{HeltonSpitkovsky2012} for further algebraic context of the numerical
range.
\par
We will discuss boundary generating curves separately for unitarily
reducible and irreducible 3-by-3 matrices. 
We will say that $\alpha\in\bC$ is a {\em normal eigenvalue} of $A$ 
is there exists a non-zero vector $x$ such that $Ax=\alpha x$,
$A^*x=\overline{\alpha}x$. The 3-by-3 matrix $A$ has a normal eigenvalue
if and only if $A$ is unitarily reducible. Every normal eigenvalue
contributes as a point to the boundary generating curve, see \S7 of
\cite{Kippenhahn1951}.
\begin{rem}\label{rem:Kippenhahn-red}
A complete list of boundary generating curves for a unitarily
reducible matrix $A\in M_3$ is
\begin{enumerate}
\item three normal eigenvalues (not necessarily distinct),
\item one normal eigenvalue and one ellipse.
\end{enumerate}
The corresponding shapes of $W(A)$ are (1) a point, a segment, or a
triangle and (2) an ellipse or the convex hull of an ellipse
and a point outside the ellipse.
\end{rem}
\begin{rem}\label{rem:Kippenhahn-irr}
A complete list of boundary generating curves for a unitarily
irreducible matrix $A\in M_3$ is
\begin{enumerate}
\item an ellipse and a point inside the ellipse,
\item a degree four curve with a double tangent,
\item a degree six curve consisting of two nested parts
one inside another, the outer part having an `ovular' shape.
\end{enumerate}
\end{rem}
\par
Following \cite{RodmanSpitkovsky2005} we denote the sets of
irreducible matrices (1) by $\cE_3$, their numerical ranges being
ellipses, and the set of irreducible matrices (2) by $\cF_3$, their
numerical ranges having a flat boundary portion. Further, we denote
the set of irreducible matrices (3) by $\cO_3$ and the set of
reducible 3-by-3 matrices by $\cR_3$. This gives a disjoint
union
\begin{equation}\label{eq:partition}
M_3=\cR_3\cup\cE_3\cup\cF_3\cup\cO_3.
\end{equation}
\par
We begin with some observations following from already published
results though not necessarily explicitly stated there.
\begin{rem}\label{rem:obs}~
\begin{enumerate}
\item The set of reducible matrices $\cR_3$ is closed and nowhere
dense in $M_3$. In fact, the analogous statement holds for all
matrix sizes $d\in\bN$, see Problem~8 in \cite{Halmos1970}.
\item
The set $\cE_3$ is closed relative to the set of unitarily
irreducible matrices, and nowhere dense. The closedness follows from
Thms.~2.3, 2.4 of \cite{Keeler-etal1997} which provide a constructive
criterion for a 3-by-3 matrix to have an elliptical numerical range.
Indeed, these conditions are in form of equalities which persist under
taking limits, while staying within the set of non-normal matrices.
This implies further that the set of 3-by-3 matrices who have an
elliptical numerical range is closed relative to the set of
non-normal matrices. The fact that $\cE_3$ is nowhere dense
was derived from the same criterion in \cite{Rault-etal2013}, see
Prop.~3.1 there.
\item
The set $\cF_3$ is closed relative to the set of unitarily
irreducible matrices and nowhere dense. The closedness follows from
Prop.~3.2 of \cite{Keeler-etal1997}, according to which a unitarily
irreducible matrix $A$ belongs to $\cF_3$ if and only if
$u\Re(A)+v\Im(A)+w\id$ has rank one for some real $u,v,w$. It remains
to invoke the lower semi-continuity of the rank
function. The fact that $\cF_3$ is nowhere dense can be seen from
an alternative description of this class, provided by
\cite[Theorem~1.2]{RodmanSpitkovsky2005}.
\item
From \eqref{eq:partition} and already established (1)--(3) it directly
follows that the set $\cO_3$ is open and dense in $M_3$.
\end{enumerate}
\end{rem}
\par
The statements (2) and (3) can be refined using Kippenhahn's
classification in Remarks~\ref{rem:Kippenhahn-red}
and~\ref{rem:Kippenhahn-irr} and the property that the numerical
ranges of a converging sequence of matrices converge to the
numerical range of the limit matrix (Rem.~\ref{rem:Hcnr}).
\begin{lem}\label{lem:refined-obs}~
\begin{enumerate}
\item
If $A\in M_3\setminus\cE_3$ is the limit of a sequence in $\cE_3$, then $A$ is reducible and $W(A)$ is an ellipse, a segment, or a point.
\item
If $A\in M_3\setminus\cF_3$ is the limit of a sequence in $\cF_3$, then $A$ is reducible. If $W(A)$ is an ellipse then the normal
eigenvalue of $A$ lies on the boundary of $W(A)$.
\end{enumerate}
\end{lem}
{\em Proof:}
(1). Rem.~\ref{rem:obs}(2) shows that $A$ is reducible and that $W(A)$ is
an ellipse if $A$ is not normal.
The numerical range $W(A)$ cannot be a triangle because the
limit of a sequence of ellipses is an ellipse, a segment, or a point
(Rem.~\ref{rem:limit-of-ellipses}).
\par
(2). Rem.~\ref{rem:obs}(3) shows that $A$ is reducible.
The Blaschke selection theorem, Thm.~18.6 of \cite{Schneider2014}, shows
that there is a sequence $(A_i)_{i\in\bN}\subset\cF_3$ converging to $A$
with flat boundary portions $s_i$ of $W(A_i)$, $i\in\bN$, such that
$(s_i)_{i\in\bN}$ converges in the Hausdorff distance. If $W(A)$ is an 
ellipse, then Rem.~\ref{rem:Hcnr}(2) shows that the limit is an extreme 
point $\alpha$ of $W(A)$ and Lemma~\ref{lem:birth-implies-multiple} 
proves that $\alpha$ is multiply generated. Now Thm.~3.2 of
\cite{Leake-etal-OAM2014} proves that $\alpha$ is an eigenvalue of $A$
which is in fact a normal eigenvalue, see Theorem 1.6.6 
of~\cite{HornJohnson1991}.
\hspace*{\fill}$\square$\\
%
%
%
%
%
\section{Reducible 3-by-3 matrices}
\label{sec:reducible}
\par
In this section we describe the intersection of the (norm) closures
$\overline{\cE_3}$ of $\cE_3$ and $\overline{\cF_3}$ of $\cF_3$
with the set of reducible matrices.
\par
To simplify possible limit points of sequences we consider the
equivalence of matrices $A\in M_d$ modulo
\begin{equation}\label{eq:condition}
\parbox{10cm}{\begin{itemize}
\item unitary similarity,
\item substitution of $A$ by a matrix in $M_d$ while preserving the
linear span of $\id,\Re(A),\Im(A)$.
\end{itemize}}
\end{equation}
The numerical range $W(A)$ and the boundary generating curve
$S_A^\wedge(\bR)$ are invariant under the action of the unitary group.
The substitutions are realized by the group of invertible affine
transformations of the real plane, whose action commutes with the
operators $A\mapsto W(A)$ and $A\mapsto S_A^\wedge(\bR)$, see
\S2.4 and~\S4.18 of \cite{Kippenhahn1951}. Hence, every equivalence class
\eqref{eq:condition} is a subset of an equivalence class in Kippenhahn's
classification provided in Remarks~\ref{rem:Kippenhahn-red}
and~\ref{rem:Kippenhahn-irr}. In particular, the blocks $\cR_3$, $\cE_3$,
$\cF_3$, and $\cO_3$ in the partition \eqref{eq:partition} of $M_d$ are
invariants of \eqref{eq:condition}. Since the action of a fixed unitary
(respectively, affine transformation) is a homeomorphism of $M_d$, the
closures $\overline{\cE_3}$ and $\overline{\cF_3}$ are invariants, too.
\begin{lem}\label{lem:red-representation}
Any matrix $A\in M_3$ such that $\Re(A)$ and $\Im(A)$
commute is equivalent modulo \eqref{eq:condition} to the
diagonal matrix
\begin{equation}\label{eq:red-comm}
0, \qquad
\diag[0,\lambda,1] \mbox{ for } 0\leq\lambda\leq\tfrac{1}{2}, \quad \mbox{or} \quad
\diag[0,1,\ii].
\end{equation}
Any reducible matrix $A\in M_3$ such that $\Re(A)$ and $\Im(A)$
do not commute is equivalent modulo \eqref{eq:condition} to
\begin{equation}\label{eq:red-nc}
\left[\begin{smallmatrix}
0 & 2 & 0\\
0 & 0 & 0\\
0 & 0 & a
\end{smallmatrix}\right],
\qquad a\geq 0.
\end{equation}
No two of the matrices in \eqref{eq:red-comm} or \eqref{eq:red-nc}
are equivalent modulo \eqref{eq:condition}.
\end{lem}
{\em Proof:}
We shall simplify $A$ using transformations \eqref{eq:condition}.
Since  $A$ is unitarily reducible we can assume that its real part
is $X\oplus a$ and its imaginary part $Y\oplus b$ where $X,Y$ are
self-adjoint 2-by-2 matrices and $a,b$ reals. Using an affine
transformation we can assume that $\tr(X)=0$, $\tr(Y)=0$, $a\geq 0$,
and $b=0$. Thereby, the real and imaginary parts of the original
matrix $A$ commute if and only if $X$ and $Y$ commute.
\par
If $X$ and $Y$ commute then, using unitary similarity, we assume
that $X$ and $Y$ are scalar multiples of $\diag[1,-1]$. The three
cases $A=0$, ($A\neq 0$, $Y=0$), and ($A\neq 0$, $Y\neq 0$)
lead {\em via} affine transformations to the cases of
\eqref{eq:red-comm} in the same order.
\par
If $X$ and $Y$ do not commute then by adding
scalar multiples of $Y\oplus 0$ to $X\oplus a$ we assume that $X$
and $Y$ are orthogonal with respect to the Hilbert-Schmidt inner
product $M_2\times M_2\to\bC$, $(x,y)\mapsto\tr(x^*y)$. Normalizing
and applying an $SU(2)=SO(3)$ rotation we obtain a matrix of the
form \eqref{eq:red-nc}.
\par
To see that no two of the matrices described in \eqref{eq:red-comm}
and \eqref{eq:red-nc} are equivalent modulo \eqref{eq:condition}
recall that the boundary generating curves of two equivalent matrices
are affinely isomorphic. The boundary generating curves of the
described matrices are one, two, or three points on a line, the
vertices of a triangle, or the union of a circle and a point. This
makes sure that no two of them are affinely isomorphic.
\hspace*{\fill}$\square$\\
\par
We shall compute the intersections of $\overline{\cE_3}$ with the set
of reducible matrices.
It is shown in \cite{RodmanSpitkovsky2005}, Sec.~2, that the
boundary generating curve of the matrix
\[
A':=
\left[\begin{smallmatrix}
a & x & y\\
0 & b & z\\
0 & 0 & c
\end{smallmatrix}\right],
\qquad a,b,c,x,y,z\in\bC
\]
consists of an ellipse and a point if and only if
$d:=|x|^2+|y|^2+|z|^2>0$ and the number
\begin{equation}\label{eq:elliptical-range-1}
\lambda:=(c|x|^2+b|y|^2+a|z|^2-x\bar yz)/d
\end{equation}
coincides with one of the eigenvalues $a,b,c$ of $A'$. In this case
the two other eigenvalues of $A'$ are the foci of an ellipse with
minor axis of length $\sqrt{d}$, and the boundary generating curve
of $A'$ is the union of this ellipse and of $\lambda$.
Moreover, the eigenvalue $\lambda$ is a normal eigenvalue
if and only if $A$ is reducible.
\par
Observe that the boundary generating curve of the following
matrix is the union of the unit circle and a point lying on
or inside the unit circle. So the numerical range is the
unit disk.
\begin{lem}\label{lem:e3-to-disks}
If $a\in[0,1]$ then the matrix
$\left[\begin{smallmatrix}
 0 & 2 & 0\\
 0 & 0 & 0\\
 0 & 0 & a
\end{smallmatrix}\right]$
lies in $\overline{\cE_3}$.
\end{lem}
{\em Proof:}
If $a>0$ then the matrix in the above statement equals
$a\cdot A(\tfrac{2}{a})$ where
\[
A(\gamma):=
\left[\begin{smallmatrix}
 0 & \gamma & 0\\
 0 & 0 & 0\\
 0 & 0 & 1
\end{smallmatrix}\right],
\qquad \gamma\geq 2.
\]
Hence, it suffices to prove $A(\gamma)\in\overline{\cE_3}$ for
real $\gamma\geq 2$. The case $a=0$ would follow by taking the
limit $a\to 0$. We define
\[
M(\alpha,\beta) :=
\left[\begin{smallmatrix}
 \alpha & (1 - \alpha)(1 + \beta^2)/\beta & \alpha\\
 0 & \alpha & -\alpha \beta\\
 0 & 0 & 1
\end{smallmatrix}\right],
\qquad
\alpha\in\bR, \beta>0.
\]
Choose $\beta>0$ such that $\gamma=(1 + \beta^2)/\beta$. Then the matrix
$A(\gamma)$ is the limit $\alpha\to 0$ of $M(\alpha,\beta)$. It suffices
to prove $M(\alpha,\beta)\in\cE_3$ for $\beta>0$ and for
$\alpha$ in a neighborhood of zero. The discussion in the paragraph of
\eqref{eq:elliptical-range-1} shows that for $\alpha>0$ the numerical
range of $M(\alpha,\beta)$ is a disk centered at $\alpha$.
\par
Let us prove that $M(\alpha,\beta)$ is unitarily irreducible for all 
$\beta>0$ and $\alpha\not\in\{0,1\}$. For $\alpha$ different from one, 
the latter is a simple eigenvalue (for $M$, as well as for $M^*$), and 
$e_3=[0,0,1]$, is the respective eigenvector for $M^*$. By 
contradiction, if $M$ were unitarily reducible, then $e_3$ would lie 
in its reducing subspace, say $L$. Applying $M$ to $e_3$ and assuming
$\alpha\neq 0$ the vector $[1,-\beta,0]$ also lies in $L$, and by 
applying $M$ again the vector $e_1=[1,0,0]$ is in $L$. So, unless 
$\beta$ is zero, $L$ is $3$-dimensional, proving that $M$ is unitarily 
irreducible.
\hspace*{\fill}$\square$\\
\begin{lem}\label{lem:e3-to-segments}
If $0\leq\lambda\leq\tfrac{1}{2}$ then the matrix $\diag[0,\lambda,1]$
lies in $\overline{\cE_3}$.
\end{lem}
{\em Proof:}
It suffices to show that for all $\lambda\in[0,\tfrac{1}{2}]$ the diagonal
matrix $\diag[0,\lambda,1]$ is a limit of matrices obtainable from
$A(a):=\left[\begin{smallmatrix}
 0 & 2 & 0\\
 0 & 0 & 0\\
 0 & 0 & a
\end{smallmatrix}\right]$, $a\in[0,1]$, considered in
Lemma~\ref{lem:e3-to-disks}, {\em via} transformations \eqref{eq:condition}.
To this end, introduce the matrix resulting from the affine transformation
$(x,y)\mapsto((1-x)/2,\epsilon y)$, $\epsilon>0$, applied to
$U^*A(1-2\lambda)U$, where
$U:=\tfrac{1}{\sqrt{2}}\left[\begin{smallmatrix}
 1 & 0 & 1\\
 1 & 0 & -1\\
 0 & \sqrt{2} & 0
\end{smallmatrix}\right]$ is unitary. A straightforward computation shows
that it equals
$\left[\begin{smallmatrix}
 0 & 0 & -\epsilon\\
 0 & \lambda & 0\\
 \epsilon & 0 & 1
\end{smallmatrix}\right]$
which converges for $\epsilon\to 0$ to $\diag[0,\lambda,1]$.
\hspace*{\fill}$\square$\\
\par
Lemmas~\ref{lem:refined-obs}(1),
\ref{lem:red-representation},
\ref{lem:e3-to-disks},
and~\ref{lem:e3-to-segments} show the following.
\begin{lem}\label{lem:E3cl}
A reducible 3-by-3 matrix lies in $\overline{\cE_3}$ if and
only if its numerical range is an ellipse, a segment, or a point.
\end{lem}
\par
We shall compute the intersection of $\overline{\cF_3}$ with the set
of reducible matrices.
\begin{lem}\label{lem:F3cl}
A reducible 3-by-3 matrix $A$ lies in $\overline{\cF_3}$ if and only if
$W(A)$ is not an ellipse having the normal eigenvalue of $A$ in the
interior.
\end{lem}
{\em Proof:}
Lemma~\ref{lem:refined-obs}(2) excludes reducible matrices $A$
from $\overline{\cF_3}$ when $W(A)$ is an ellipse with the normal
eigenvalue of $A$ in the interior. Let us show that all other
reducible matrices lie in $\overline{\cF_3}$.
\par
We can assume that $A$ is a matrix listed in
Lemma~\ref{lem:red-representation}. Let $A$ be of the form
\eqref{eq:red-comm}. Clearly $0\in\overline{\cF_3}$. Otherwise,
if $A\neq 0$, then $A=\diag[0,\lambda,1]$ for $\lambda\in[0,\tfrac{1}{2}]$
or $A=\diag[0,1,\ii]$. In both cases we define for $\epsilon>0$
a matrix $M(\epsilon)_{j,k}:=A_{j,k}+\ii\epsilon$, $j,k=1,2,3$. In
the first case $\Im(M(\epsilon))$ and in the second case
$\Re(M(\epsilon))$ has a multiple eigenvalue. If $\epsilon>0$
then these matrices are irreducible because their real and imaginary
parts have no common eigenvectors. This proves $A\in\overline{\cF_3}$.
\par
Let $A$ be of the form \eqref{eq:red-nc} with $a\geq 1$.
It suffices to consider $a>1$, the case $a=1$ would follow by
taking the limit. We apply transformations
\eqref{eq:condition}, namely the affine transformation
$\bR^2\to\bR^2$, $[x,y]\mapsto[x+y\sqrt{a^2-1},y]$ followed by
a unitary transformation $A\mapsto u^*Au$ with unitary
\[
u:=\tfrac{1}{a\sqrt{2}}
\left[\begin{smallmatrix}
 0 & 1-\ii \sqrt{a^2-1} & 1-\ii \sqrt{a^2-1} \\
 0 & a & -a \\
 a\sqrt{2} & 0 & 0
\end{smallmatrix}\right].
\]
By means of these transformations we may assume that
\[
\Re(A)=\diag[a,a,-a]
\quad\mbox{and}\quad
\Im(A)=\tfrac{1}{a}
\left[\begin{smallmatrix}
 0 & 0 & 0 \\
 0 & \sqrt{a^2-1} & \ii \\
 0 & -\ii & -\sqrt{a^2-1}
\end{smallmatrix}\right].
\]
Now it is obvious that for all $\epsilon>0$ the matrix
\[
\diag[a,a,-a] + \tfrac{\ii}{a}
\left[\begin{smallmatrix}
 0 & 0 & \epsilon \\
 0 & \sqrt{a^2-1} & \ii \\
 \epsilon & -\ii & -\sqrt{a^2-1}
\end{smallmatrix}\right]
\]
is unitarily irreducible and Rem.~\ref{rem:obs}(3) shows that
it lies in $\cF_3$. Since the matrix converges to $A$ for
$\epsilon\to 0$, the proof is complete.
\hspace*{\fill}$\square$\\
\par
Lemmas~\ref{lem:refined-obs}, \ref{lem:E3cl} and \ref{lem:F3cl}
combined yield the following.
\begin{thm}
The set $\overline{\cE_3}\cap\overline{\cF_3}$ is the
subset of all reducible matrices $A\in\cR_3$ where $W(A)$
is a point, a segment, or an ellipse with the normal
eigenvalue of $A$ on the boundary.
\end{thm}
%

%
%
%
%

\begin{thebibliography}{10}

\bibitem{Au-YeungPoon1979} Y.\,H.~Au-Yeung, Y.\,T.~Poon (1979)
{\em A remark on the convexity and positive definiteness
concerning Hermitian matrices},
Southeast Asian Bull.\ Math.\ 3 85--92

\bibitem{BerberianOrland1967} S.\,K.~Berberian, G.\,H.~Orland (1967)
{\em On the closure of the numerical range of an operator},
Proceedings of the American Mathematical Society 18(3) 499--503

\bibitem{BonnesenFenchel1987} T.~Bonnesen, W.~Fenchel (1987)
{\em Theory of Convex Bodies}, BCS Associates, Moscow, Idaho USA

\bibitem{BrownSpitkovsky2004} E.\,S.~Brown, I.\,M.~Spitkovsky (2004)
{\em On flat portions on the boundary of the numerical range},
Linear Algebra and its Applications 390 75--109

\bibitem{Chen-etal2015} J.~Chen, Z.~Ji, C.-K.~Li, Y.-T.~Poon, Y.~Shen,
N.~Yu, B.~Zeng, D.~Zhou (2015)
{\em Discontinuity of maximum entropy inference and quantum phase transitions},
New Journal of Physics 17(8) 083019

\bibitem{ChienNakazato2010} M.-T.~Chien, H.~Nakazato (2010)
{\em Joint numerical range and its generating hypersurface},
Linear Algebra and its Applications 432(1) 173--179

\bibitem{Corey-etal2013} D.~Corey, C.\,R.~Johnson, R.~Kirk,
B.~Lins, I.~Spitkovsky (2013)
{\em Continuity properties of vectors realizing points in the
classical field of values},
Linear and Multilinear Algebra 61 1329--1338

\bibitem{CsiszarMatus2005} I.~Csisz\'ar, F.~Mat\'u\v s (2005)
{\em Closures of exponential families},
The Annals of Probability 33(2) 582--600

\bibitem{Eldred-etal2012} J.~Eldred, L.~Rodman, I.~Spitkovsky (2012)
{\em Numerical ranges of companion matrices: flat portions on the
boundary}, Linear and Multilinear Algebra 60 1295--1311

\bibitem{Fischer2001} G.~Fischer (2001)
{\em Plane Algebraic Curves}, Providence, Rhode Island: AMS

\bibitem{GallaySerre2012} T.~Gallay, D.~Serre (2012)
{\em Numerical measure of a complex matrix},
Communications on Pure and Applied Mathematics 65(3) 287--336

\bibitem{Gelfand-etal1994} I.\,M.~Gelfand, M.\,M.~Kapranov,
A.\,V.~Zelevinsky (1994)
{\em Discriminants, Resultants, and Multidimensional Determinants},
Boston, MA: Birkh\"auser Boston

\bibitem{Gruenbaum2003} B.~Gr\"unbaum (2003)
{\em Convex Polytopes}, 2nd Edition, New York: Springer

\bibitem{Halmos1970} P.\,R.~Halmos (1970)
{\em Ten problems in Hilbert space},
Bull.\ Amer.\ Math.\ Soc.\ 76(5) 887--933

\bibitem{Hausdorff1919} F.~Hausdorff (1919)
{\em Der Wertvorrat einer Bilinearform},
Math.~Z.\ 3(1) 314--316

\bibitem{HeltonSpitkovsky2012} J.\,W.~Helton, I.\,M.~Spitkovsky (2012)
{\em The possible shapes of numerical ranges},
Operators and Matrices 6 607--611

\bibitem{HornJohnson1991} R.\,A.~Horn, C.\,R. Johnson (1991)
{\em Topics in Matrix Analysis}, 10.\ printing,
Cambridge Univ.\ Press

\bibitem{Jaynes1957} E.\,T.~Jaynes (1957)
{\em Information theory and statistical mechanics},
Phys Rev 106 620--630 and 108 171--190

\bibitem{JoswigStraub1998} M.~Joswig, B.~Straub (1998)
{\em On the numerical range map},
Journal of the Australian Mathematical Society 65 267--283

\bibitem{Kato-etal2016} K.~Kato, F.~Furrer, M.~Murao (2016)
{\em Information-theoretical analysis of topological entanglement 
entropy and multipartite correlations},
Physical Review A 93 022317

\bibitem{Keeler-etal1997} D.\,S.~Keeler, L.~Rodman,
I.\,M.~Spitkovsky (1997)
{\em The numerical range of $3\times 3$ matrices},
Lin Alg Appl 252 115--139

\bibitem{Kippenhahn1951} R.~Kippenhahn (1951)
{\em \"Uber den Wertevorrat einer Matrix},
Math Nachr 6 193--228

\bibitem{Leake-etal-OAM2014}
T.~Leake, B.~Lins, I.\,M.~Spitkovsky (2014)
{\em Pre-images of boundary points of the numerical range},
Operators and Matrices 8 699--724

\bibitem{Leake-etal-LMA2014}
T.~Leake, B.~Lins, I.\,M.~Spitkovsky (2014)
{\em Inverse continuity on the boundary of the numerical range},
Linear and Multilinear Algebra 62 1335--1345

\bibitem{LiPoon2000} C.-K.~Li, Y.-T.~Poon (2000)
{\em Convexity of the joint numerical range},
SIAM J Matrix Anal A 21(2) 668--678

\bibitem{Rault-etal2013} P.\,X.~Rault, T.~Sendova, I.\,M.~Spitkovsky (2013)
{\em 3-by-3 matrices with elliptical numerical range revisited},
Electronic Journal of Linear Algebra 26 158--167

\bibitem{Rellich1954} F.~Rellich (1954)
{\em Perturbation Theory of Eigenvalue Problems},
Research in the Field of Perturbation Theory and Linear Operators,
Technical Report No.~1,
Courant Institute of Mathematical Sciences, New York University

\bibitem{RodmanSpitkovsky2005} L.~Rodman, I.\,M.~Spitkovsky (2005)
{\em $3\times 3$ matrices with a flat portion on the boundary of the
numerical range},
Linear Algebra and its Applications 397 193--207

\bibitem{Rodman-etal2016} L.~Rodman, I.\,M.~Spitkovsky,
A.~Szko{\l}a, S.~Weis (2016)
{\em Continuity of the maximum-entropy inference:\
Convex geometry and numerical ranges approach},
Journal of Mathematical Physics 57, 015204

\bibitem{Schneider2014} R.~Schneider (2014)
{\em Convex Bodies:\ The Brunn-Minkowski Theory}, 2nd Edition,
Cambridge University Press

\bibitem{Toeplitz1918} O.~Toeplitz (1918)
{\em Das algebraische Analogon zu einem Satze von Fej\'er},
Math.~Z.\ 2(1--2) 187--197

\bibitem{Weis-support} S.~Weis (2011)
{\em Quantum convex support},
Linear Algebra and its Applications
435(12) 3168--3188; correction (2012) ibid.\ 436 xvi

\bibitem{Weis-cones} S.~Weis (2012)
{\em A note on touching cones and faces},
Journal of Convex Analysis 19 323--353

\bibitem{Weis-topo} S.~Weis (2014)
{\em Information topologies on non-commutative state spaces},
Journal of Convex Analysis 21 339--399

\bibitem{Weis-cont} S.~Weis (2014)
{\em Continuity of the maximum-entropy inference},
Communications in Mathematical Physics 330(3) 1263--1292

\bibitem{Weis-strong} S.~Weis (in print)
{\em Maximum-entropy inference and inverse continuity of
the numerical range}, Reports on Mathematical Physics 

\bibitem{WeisKnauf2012} S.~Weis, A.~Knauf (2012)
{\em Entropy distance:\ New quantum phenomena},
J Math Phys 53(10) 102206

\end{thebibliography}
\bibliographystyle{plain}

\vspace{.5in}

\begin{center}
\begin{tabular}{l}
%
%
%
%
%
Ilya M.\ Spitkovsky \\
Division of Science and Mathematics \\
New York University Abu Dhabi \\
Saadiyat Island, P.O. Box 129188 \\
Abu Dhabi, UAE \\
e-mail: {\tt ims2@nyu.edu}\\
~\\
and\\
~\\
Department of Mathematics \\
College of William and Mary \\
P.\, O.~Box 8795 \\
Williamsburg, VA 23187-8795 \\
e-mail: {\tt ilya@math.wm.edu}\\
~\\
Stephan Weis\\ 
Departamento de Matemática\\
Universidade Estadual de Campinas\\			
Rua Sérgio Buarque de Holanda, 651\\
Campinas-SP, 13.083-859\\ 
Brazil\\
e-mail: {\tt maths@stephan-weis.info}
\end{tabular}
\end{center}
\end{document}